\documentclass{amsart}

\usepackage[utf8]{inputenc}
\usepackage[T1]{fontenc}
\usepackage{natbib}

\usepackage{amsmath}
\usepackage{amssymb}
\usepackage{amsthm}
\usepackage{amsfonts}
\usepackage{mathrsfs}
\usepackage{xspace}
\usepackage{tikz}
\usepackage{fixmath}
\usepackage{mathtools}
\usepackage{array}

\usepackage{fixltx2e}
\usepackage[expansion=false,final]{microtype}

\newtheorem{thm}{Theorem}[section]
\newtheorem{lem}[thm]{Lemma}
\newtheorem{defn}[thm]{Definition}
\newtheorem{prop}[thm]{Proposition}
\newtheorem{cor}[thm]{Corollary}

\DeclareMathOperator{\Ext}{Ext}

\DeclareMathOperator{\ch}{ch}
\DeclareMathOperator{\todd}{Td}
\DeclareMathOperator{\rk}{rk}
\DeclareMathOperator{\Hom}{Hom}
\DeclareMathOperator{\chern}{ch}
\DeclareMathOperator{\Pic}{Pic}
\DeclareMathOperator{\Arg}{arg}
\DeclareMathOperator{\tor}{tor}
\DeclareMathOperator{\re}{Re}
\DeclareMathOperator{\im}{Im}

\newcommand{\Z}{\mathbb{Z}\xspace}

\newcommand{\cL}{\mathcal{L}\xspace}
\newcommand{\cM}{\mathcal{M}\xspace}

\newcommand{\cI}{\mathcal{I}\xspace}
\newcommand{\cH}{\mathcal{H}\xspace}
\newcommand{\R}{\mathbb{R}\xspace}
\newcommand{\C}{\mathbb{C}\xspace}

\newcommand{\E}{\mathbb{E}\xspace}
\newcommand{\cP}{\mathcal{P}\xspace}

\newcommand{\cF}{\mathcal{F}\xspace}
\newcommand{\cA}{\mathcal{A}\xspace}
\newcommand{\cG}{\mathcal{G}\xspace}
\newcommand{\cQ}{\mathcal{Q}\xspace}

\newcommand{\bP}{\mathbb{P}\xspace}
\newcommand{\sC}{\mathscr{C}\xspace}
\newcommand{\cO}{\mathcal{O}\xspace}

\title{Nef and Effective Cones of the Moduli Spaces of Torsion Sheaves on \(\mathbb{P}^2\)}
\author{Matthew Woolf}

\begin{document}

\bibliographystyle{chicago}

\begin{abstract}
In this paper, we study the divisor theory of the Simpson moduli space of semistable sheaves of dimension 1 on \(\mathbb{P}^2\). We prove that these spaces are all Mori dream spaces, and calculate their nef cones. We also study the effective cones of these spaces for most choices of numerical invariants. As a consequence, we work out precisely when two such spaces are isomorphic.
\end{abstract}

\maketitle

\tableofcontents

\section{Introduction}

The moduli space of torsion sheaves on \(\bP^2\) is a space parametrizing semistable sheaves with Hilbert polynomial \(\chi(\cF(m))=\mu m+\chi\) on $\bP^2$ (technically, their S-equivalence classes). In this paper, we will call this space $N(\mu,\chi)$. This space is a compactification of the relative Picard scheme for line bundles of Euler characteristic $\chi$ on a smooth plane curve of degree $\mu$. In this paper, we will study the divisor theory of these moduli spaces, specifically, their pseudoeffective and nef cones.

When $\mu \geq 3$, these spaces all have Picard number 2, so their divisor theory is relatively simple. Moreover, we will show that they are Mori dream spaces, so even the pathologies which can occur for a space of Picard number 2 are not a problem for these moduli spaces.

On these moduli spaces, there is always a natural divisor class which lies on the edge of both the nef and effective cones. Specifically, there is a map from $N(\mu,\chi)$ to the projective space of all curves of degree $\mu$. We can then pull back $\cO(1)$ from this projective space to get such a divisor class. Much of this paper is devoted to finding the other edge of both the nef and the pseudoeffective cones.

In doing so, we find something similar to what is known for the Hilbert scheme. That is, the effective cones is related to the algebra of these sheaves, while the nef cone is more closely related to the geometry of these sheaves, when they are thought of as line bundles on plane curves. To be more specific, in the case of the effective cone, we will prove that for many of these moduli spaces, there is a vector bundle $E$ such that the locus of $\cF \in N(\mu,\chi)$ such that $h^0(\E \otimes \cF)$ is a divisor, and it is this divisor which lies on the other edge of the effective cone. The divisor on the other edge of the nef cone is the pullback of an ample divsior by a regular map, and we will describe the fibers of this map in terms of points on curves and the line bundles they define.

For the case of $\chi=1$, Jinwon Choi and Kiryong Chung have a much more geometric description of the divisor lying on the other edge of the effective cone in \citep{choichung}.

This paper belongs to the growing genre studying the explicit birational geometry of moduli spaces. It draws much of its inspiration from \citep{abch}, which suggests that the birational models of the Hilbert scheme of points on $\bP^2$ can be realized as moduli spaces of Bridgeland stable objects. This paper aims to expand upon that one in two ways. First, it proves a number of facts about the moduli space of torsion sheaves which are hard to relate to the language of stability conditions, for example, the proof that the moduli space is a Mori dream space, or the calculation of the effective cone. Second, it extends a number of the general results which in \citep{abch} are proved for torsion-free sheaves to the case of torsion sheaves, and uses these results to carry out certain explicit calculations, for example the calculation of the nef cone.

The spaces $N(\mu,\chi)$ and $N(\mu,\chi')$ are known to be isomorphic if $\chi \cong \pm \chi' \pmod \mu$. Intuitively, these isomorphisms come from taking the tensor product of a sheaf with some multiple of $\cO(1)$, and from taking a line bundle supported on some plane curve to its dual. If $\mu \neq \mu'$, then the spaces $N(\mu,\chi)$ and $N(\mu',\chi')$ do not have the same dimension, so they cannot be isomorphic. Using the calculation of the nef cone, we will show that if $\chi \not \cong \pm \chi' \pmod \mu$ and $\mu>2$, then the dimension of the exceptional locus of a certain canonically defined map distinguishes between $N(\mu,\chi)$ and $N(\mu,\chi')$.

\textbf{Acknowledgements.} I would like to thank Dawei Chen, Izzet Coskun, Joe Harris, and Jack Huizenga for the great help they've provided me on this project. I would also like to thank Arend Bayer, Aaron Bertram, Tom Bridgeland, Gabriel Bujokas, and Emanuele Macri for very useful conversations.

\section{Preliminaries on the Moduli Space of Torsion Sheaves}

\begin{defn}
A coherent sheaf on $\bP^2$ is called pure of dimension 1 if its support is one-dimensional, and the same is true of any nonzero subsheaf.

A pure one-dimensional sheaf $\cF$ on $\bP^2$ is called semistable if for all nonzero proper subsheaves $\cG \subset \cF$, we have \[ \frac{\chi(\cG)}{\ch_1(\cG)} \leq \frac{\chi(\cF)}{\ch_1(\cF)} \] and it is called stable if moreover this inequality is always strict. This notion of stability will sometimes be called Simpson (semi)stability. (This agrees with the usual definition of Gieseker-Simpson stability applied to pure one-dimensional sheaves.)
\end{defn}

In this section, we will recall a number of facts about such sheaves, most of which can be found in le Potier's paper \citep{lepotier}, except where otherwise noted.

Given any semistable sheaf $\cF$, there is always a filtration \[ 0=\cF_0 \subset \cF_1 \subset \cdots \subset \cF_n=\cF \] called the Jordan-H\"older filtration, such that the subquotients $\cF_i/\cF_{i-1}$ are stable. This filtration is not unique, but the subquotients are determined up to reordering. Two semistable sheaves are called S-equivalent if the stable sheaves appearing as subquotients in their Jordan-H\"older filtrations are the same. We note that this filtration is only interesting if the sheaf is strictly semistable, i.e.~semistable but not stable, so each stable sheaf is the only object in its S-equivalence class.

There is a moduli space, which we will call $N(\mu,\chi)$, parametrizing S-equivalence classes of pure one-dimensional sheaves which are semistable and have Hilbert polynomial $\chi(\cF(m))=\mu m+\chi$. We will call this space the moduli space of one-dimensional sheaves, or the moduli space of torsion sheaves. In a slight abuse of notation, we will write $\cF \in N(\mu,\chi)$ to mean that $\cF$ is a semistable pure one-dimensional sheaf with Hilbert polynomial $\chi(\cF(m))=\mu M+\chi$. If $c \in K(\bP^2)$ is the class of a pure one-dimensional sheaf, then $N(c)$ will be the moduli space of semistable sheaves with class $c$, which will be isomorphic to some $N(\mu,\chi)$.

The spaces $N(\mu,\chi)$ are irreducible and factorial of dimension $\mu^2+1$. A dimension count then shows that the generic such sheaf is the push forward of a line bundle on a smooth plane curve. The spaces $N(\mu,\chi)$ are smooth away from the locus of strictly semistable sheaves. The space $N(\mu,\chi)$ is isomorphic to $N(\mu,\chi+\mu)$, with the isomorphism being given by the map $\cF \mapsto \cF \otimes \cO_{\bP^2}(1)$. We also know that the space $N(\mu,\chi)$ is isomorphic to $N(\mu,-\chi)$ by the map $\cF \mapsto \mathscr Ext^1(\cF,\omega_{\bP^2})$, as proved in \citep{maican}.

\begin{cor}\label{isomorphisms}
If $\chi \cong \pm \chi' \pmod \mu$, then $N(\mu,\chi) \cong N(\chi,\chi')$.
\end{cor}

When $\mu \geq 3$, the Picard group of $N(r,\chi)$ is a free abelian group of rank 2. We will now describe generators for this group. There is a morphism $N(\mu,\chi) \to \bP(H^0(\cO_{\bP^2}(\mu)))$ which sends a sheaf to its support. Here, the scheme-theoretic structure on the support is that given by the Fitting ideal of the sheaf. The pullback of $\cO(1)$ from $\bP(H^0(\cO_{\bP^2}(\mu)))$ to $N(\mu,\chi)$ will be denoted $\cL_0$. This is one generator of the Picard group.

Let $c$ be the class of an element of $N(r,\chi)$ in $K(\bP^2)$, the Grothendieck group of coherent sheaves on $\bP^2$. Let $a \in K(\bP^2)$. Let $\mathcal{F}$ be a sheaf on $\bP^2 \times S$ flat over $S$ such the the restriction to each fiber has class $c$, i.e.~a flat family of sheaves of class $c$ on $\bP^2$ parametrized by $S$. Let $\pi_1$ (resp.~$\pi_2$) be the projection from $\bP^2 \times S$ to $\bP^2$ (resp.~$S$). Then \[ \det(\sum_i (-1)^i R\pi_{2*}(\cF \otimes \pi_1^*(a))) \] is a line bundle on $S$.

On $K(\bP^2)$, there is a quadratic form given by \[ \langle [E],[E'] \rangle \mapsto \sum_i (-1)^i \chi(Tor^i(E,E')). \] 
Given $b \in K(\bP^2)$, let $b^\perp$ be the orthogonal complement to $b$ with respect to this quadratic form. If $a \in c^{\perp}$, then the line bundle above descends to a line bundle on $N(\mu,\chi)$. This determines a group homomorphism $\lambda: c^{\perp} \to \Pic(N(\mu,\chi))$.

\begin{defn}
The line bundle $\lambda(a)$ is called the determinant line bundle determined by $a$.
\end{defn}

As a first case, we note that $\lambda(-h^2)=\lambda_0$.

Let $\delta=\gcd(\mu,\chi)$. Let $\cL_1$ be the line bundle on $N(\mu,\chi)$ corresponding to the class \[ a=\frac{1}{\delta}((-\mu)+\chi h) \in K(\bP^2) \] where $h=[\cO_H]$. Then $\cL_0$ and $\cL_1$ generate the Picard group of $N(\mu,\chi)$. It is clear that $\cL_0$ gives one edge of both the nef and effective cones, since it is the pullback of an ample line bundle by a nonconstant regular map to a variety of smaller dimension (when $\mu \geq 3$).
 
\section{The Moduli Spaces are Mori Dream Spaces}

In this section, we will prove that the moduli spaces $N(\mu,\chi)$ are Mori dream spaces. The case when the Picard group is isomorphic to $\Z$ is trivial, so for the rest of this section, we will focus on the case $\mu \geq 3$, when the Picard group has rank two. Our first goal will be to understand the canonical class.

\begin{lem}
The canonical class $K$ of $N(\mu,\chi)$ is a negative multiple of $\cL_0$.
\end{lem}
\begin{proof}
This follows from theorem 8.3.3 of \citep{huybrechtslehn} and the preceding discussion. In fact, $K=-3\mu\cL_0$. Note that the divisor they call $\cL_1$ is in fact a multiple of what we call $\cL_0$.
\end{proof}

We will now work to understand the singularities of $N(\mu,\chi)$. The first thing to note is that as proved in sections 4.3 and 4.4 of \citep{huybrechtslehn}, for all $m$ sufficiently large, there is a $k$ such that $N(\mu,\chi)$ is the good quotient of an open subscheme $U_m$ of the scheme parametrizing quotients of $\cO(-m)^k$ with Hilbert polynomial $\mu m+\chi$. 

Write $\cH$ for $\cO(-m)^k$ and $K$ for the kernel of some surjective map $\cH \to \cF$, with $\cF$ some fixed semistable sheaf with Hilbert polynomial $\chi(\cF(n))=\mu n + \chi$. By the same discussion in \citep{huybrechtslehn}, $m$ can be picked larger than the Castelnuovo-Mumford regularity of $\cF$ for all semistable $\cF$ with this Hilbert polynomial.

\begin{lem}
The scheme $U_m$ is nonsingular for $m$ sufficiently large.
\end{lem}
\begin{proof}

By proposition 4.4.4 of \citep{sernesi}, the obstruction space to the point of $U_m$ corresponding to the quotient $\cH \to \cF$ is given by $\Ext^1(K,\cF)$, so it suffices to show that this space vanishes.

We have the short exact sequence
\[ 0 \to K \to \cH \to \cF \to 0 \]

Applying $\Hom(\cdot,\cF)$, we get a long exact sequence of $\Ext$ groups, which in particular gives us the exact sequence \[ \Ext^1(\cH,\cF) \to \Ext^1(K,\cF) \to \Ext^2(\cF,\cF) \] so it suffices to show that the two outside groups vanish.

We have $\Ext^1(\cH,\cF) \cong H^1(\cF(m))^k$, but since $m$ is larger than the regularity of $\cF$, $H^1(\cF(m))=0$. By Serre duality, $\Ext^2(\cF,\cF) \cong \Hom(\cF,\cF(-3))^*$. But $\cF(-3)$ is a pure 1-dimensional semistable sheaf with the same first Chern class as $\cF$ and smaller Euler characteristic, so by considering the image of a putative map $\cF \to \cF(-3)$, we will be able to contradict the semistability of $\cF$ unless this map is 0.
\end{proof}

By \citep{boutot}, a good quotient of a variety with rational singularities (e.g.~a nonsingular variety) also has rational singularities. This shows that $N(\mu,\chi)$ has rational (in particular, Cohen-Macaulay) singularities. By \citep{murthy}, a local ring is Gorenstein if it is Cohen-Macaulay, factorial, and a quotient of a regular local UFD. This last property will certainly be satisfied for any local ring occurring here, so we have shown that $N(\mu,\chi)$ has only Gorenstein singularities.

By theorem 11.1 of \citep{kollar} (noting that a variety with Gorenstein singularities has an invertible dualizing sheaf), rational Gorenstein singularities are canonical, so we have proved the following result.

\begin{prop}
The space $N(\mu,\chi)$ has only canonical (and hence klt, since there is no boundary divisor) singularities.
\end{prop}

By corollary 2.35 of \citep{kollarmori}, this implies that if $D$ is an effective divisor, and $\epsilon>0$ is sufficiently small, then $(N(\mu,\chi),\epsilon D)$ is a klt pair.

\begin{defn}
If $D_1,\ldots,D_n$ are Cartier divisors on a variety $X$, then we will denote by $R(X;D_1,\ldots,D_n)$ the multisection ring of the $D_i$, i.e.\[ R(X;D_1,\ldots,D_n)=\bigoplus_{(i_1,\ldots,i_n) \in \mathbb{N}^n} H^0(i_1D_1 + \cdots i_n D_n) \] Where there is no ambiguity, we will drop the $X$ from the notation.
\end{defn}

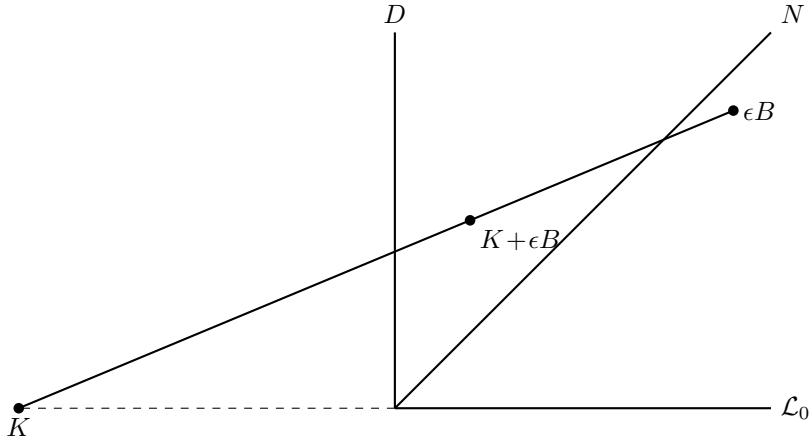
\begin{tikzpicture}[
scale=5,
axis/.style={thick},
dashed line/.style={dashed,thin},
important line/.style={thick}]
\draw[axis] (0,0) -- (1,0) node (L0line) [right] {$\cL_0$};
\draw[axis] (0,0) -- (0,1) node (effline) [above] {$D$};
\draw[axis] (0,0) -- (1,1) node (nefline) [above right] {$N$};
\draw[dashed line] (0,0) -- (-1,0) node (Canonical) [below] {$K$};
\draw[important line] (-1,0) -- (0.2,0.5) node (Boundary) [below right] {$K\!+\!\epsilon B$};
\draw[important line] (0.2,0.5) -- (0.9,0.792) node (Ample) [right] {$\epsilon B$}; 
\fill[black] (-1,0) circle (.4pt);
\fill[black] (0.2,0.5) circle (.4pt);
\fill[black] (0.9,0.792) circle (.4pt);
\end{tikzpicture}

Let $A$ be an ample divisor on $N(\mu,\chi)$. Let $D$ be a divisor on the edge of the pseudoeffective cone which does not consist of multiples of $\cL_0$. Let $N$ be a divisor on the edge of the nef cone which does not consist of multiples of $\cL_0$ By a result of Zariski (lemma 2.8 of \citep{hukeel}), since $\cL_0$ and $A$ are both semiample, $R(\cL_0,A)$ is finitely generated. Note that any effective divisor which is not a multiple of $\cL_0$ can be written as $K+\epsilon B$, with $B$ ample and $\epsilon$ arbitrarily small. By choosing $B$ so that $\epsilon$ is sufficiently small, we can ensure that $(N(\mu,\chi),\epsilon B)$ a klt pair. By corollary 1.19 of \citep{bchm}, $R(A,D)$ is finitely generated.

\begin{prop}\label{moridreamspace}
The Cox ring of $N(\mu,\chi)$ is finitely generated.
\end{prop}
\begin{proof}
We note that the Cox ring of $N(\mu,\chi)$ is isomorphic to $R(\cL_0,D)$. Inside $R(\cL_0,D)$, we have two finitely generated subrings, namely $R(\cL_0,A)$ and $R(A,D)$. I claim that the generators of these two rings together generate $R(\cL_0,D)$. It clearly suffices to show that any section of any effective divisor $D$ on $N(\mu,\chi)$ is in $R(\cL_0,A)$ or $R(A,D)$. But $D$ is either in the cone generated by $\cL_0$ and $A$ or it is in the cone generated by $A$ and $D$. In the former case, the section of $D$ lies in $R(\cL_0,A)$, and in the latter case, it is contained in $R(A,D)$.
\end{proof}

By \citep{hukeel}, we have the following corollary.

\begin{cor}
The spaces $N(\mu,\chi)$ are Mori dream spaces.
\end{cor}

\section{Effective Divisors}

In the case of $N(\mu,0)$, the line bundle $\cL_1=\lambda(-1)$ is effective -- it corresponds to the subvariety of $N(\mu,\chi)$ consisting of sheaves which have a global section (\citep{lepotier}). More generally, we have the following.

\begin{defn}
Let $\cF$ be a sheaf on $\bP^2$ and $E$ a vector bundle. We say that $E$ is cohomologically orthogonal to $\cF$ if $H^i(\cF \otimes E)=0$ for all $i$. If there is some semistable $\cF \in N(\mu,\chi)$ which is cohomologically orthogonal to $E$, we say that $E$ has the interpolation property for $N(\mu,\chi)$. We will also talk about the interpolation property for the Hilbert scheme of points on $\bP^2$, which has a completely analagous definition, since the Hilbert scheme of $n$ points on $\bP^2$ is isomorphic to a moduli space of semistable sheaves on $\bP^2$ with rank 1, first Chern class 0, and second Chern class $-n$ (see for example \citep{huizenga2}).
\end{defn}

\begin{prop}
Let $c \in K(\bP^2)$ denote the class of some semistable pure one-dimensional sheaf $\cF$. Suppose $E$ is cohomologically orthogonal to $\cF$. Let $D$ be the locus in $N(c)$ of sheaves which are not cohomologically orthogonal to $E$. Then $D$ is a divisor on $N(\mu,\chi)$, and the corresponding line bundle is the dual of the determinant bundle corresponding to $E$.
\end{prop}

The proof of this proposition follows from the same argument as in the discussion after proposition 2.10 of \citep{lepotier}, which covers the case where $E=\cO_{\bP^2}$. A calculation with the Hirzebruch-Riemann-Roch formula shows that $\chi(E,\cF)=0$, which is a necessary condition for $E$ to be cohomologically orthogonal to $\cF$, if and only if the slope of $E$ is $\frac{-\chi}{\mu}$.

In particular, since the map $\lambda:c^\perp \to \Pic(N(\mu,\chi))$ is a homomorphism, the class of $D$ is $a \cL_0+b \cL_1$, where $a$ and $b$ are such that \[ \frac{b}{\delta}(-\mu+\chi \chern\cO_H)+a h^2=-\chern(E) \] More explicitly, we get $a=-\ch_2E+\frac{\chi}{2\mu} \rk E$ and $b=\frac{\delta}{\mu} \rk(E)$.

Because of the isomorphism $N(\mu,\chi) \cong N(\mu,\chi+\mu)$, in order to calculate the effective cone, it suffices to assume that $0<\chi \leq \mu$. The isomorphism $N(\mu,\chi) \cong N(\mu,-\chi)$ lets us reduce further to the case $0 \leq \chi \leq \frac{\mu}{2}$. We will find the effective cone in these cases (when we can) by comparing the minimal resolutions of generic one-dimensional sheaves and of generic ideal sheaves of zero-dimensional subschemes. We will use the similarity between the minimal resolutions to show that there is a relation between cohomological orthogonality in the two cases for a particular type of vector bundle.

We now consider the Hilbert scheme $\cH_n$ of $n=\binom{\mu+1}{2}+\mu-\chi$ points in $\bP^2$.

\begin{thm}\label{steiner}
Let $0 \leq 2\chi \leq \mu$. If we have $E$ a vector bundle which fits in a short exact sequence \[ 0 \to k(\mu-\chi)\cO(\mu-2) \to k(2\mu-\chi)\cO(\mu-1) \to E(\mu) \to 0 \] for some positive integer $k$ such that $E(\mu)$ has the interpolation property with respect to the Hilbert scheme of $n$ points in $\bP^2$, then $E$ has the interpolation property with respect to $N(\mu,\chi)$.
\end{thm}

We note that the class of the divisor corresponding to such an $E$ is \[ D=k((\mu-\chi) \cL_0+\delta \cL_1) \]

To prove this theorem, we first need to recall the following results.

\begin{prop}
Let $\cF$ be a general element of $N(\mu,\chi)$, with $0 \leq 2\chi \leq \mu$. The minimal free resolution of $\cF$ is \[ 0 \to (\mu-\chi) \cO_{\bP^2}(-2) \to \chi \cO_{\bP^2} \oplus (\mu-2\chi)\cO_{\bP^2}(-1) \to \cF \to 0. \]
%
\end{prop}

This is proved in \citep{freiermuth}. Moreover, by openness of semistability and injectivity, if we take a general such map, the cokernel will be a semistable sheaf with the correct invariants. Note that this implies that $N(\mu,\chi)$ is unirational.

\begin{prop}
Let $\cI$ be a general element of $\cH_n$ and $\cI'=\cI(\mu)$. If $0 \leq 2\chi \leq \mu$, then the minimal free resolution of $\cI'$ is  \[ 0 \to (\mu-\chi) \cO_{\bP^2}(-2) \to (\chi+1) \cO_{\bP^2} \oplus (\mu-2\chi)\cO_{\bP^2}(-1) \to \cI' \to 0. \]
%
\end{prop}

This is proved in \citep{eisenbud}. Moreover, if we take the cokernel of a general such map, we get a (twist of) an ideal sheaf of the right numerical invariants.

We can now prove the above theorem.

\begin{proof} 
Let $\cI$ be a general ideal sheaf in $\cH_n$, and $\cI'=\cI(\mu)$. We note that $E$ has slope $\frac{-\chi}{\mu}$, which means that $\chi(E,\cF)=0$ for $\cF \in N(\mu,\chi)$. It therefore suffices to show that $H^0(\cF \otimes E) \cong H^0(\cI' \otimes E)$, since $H^0(E \otimes \cI')=0$ for general $\cI'$ by assumption.

Note that by the long exact sequence in cohomology coming from \[ 0 \to a\cO(-2) \to b\cO(-1) \to E \to 0 \] we see that $H^i(E)=0$ for all $i$. Since we're assuming that $E(\mu)$ is cohomologically orthogonal to the general ideal sheaf, we see that given a general injective map \[ (\mu-\chi) \cO_{\bP^2}(-2) \to (\mu-2\chi) \cO_{\bP^2}(-1) \oplus (\chi+1) \cO_{\bP^2} \] if we denote the cokernel of the map by $C$, then $H^i(E \otimes C)=0$ for all $i$. Using the long exact sequence in cohomology coming from the short exact sequence defining $C$, we see that this is equivalent to the maps \[(\mu-\chi) H^i(E(-2)) \to (\chi+1) H^i(E) \oplus (\mu-2\chi) H^i(E(-1))  \] being isomorphisms.

We know that $H^i(E)=0$ for all $i$ by the long exact in cohomology coming from the defining short exact sequence of $E$, the cohomological orthogonality of $E(\mu)$ to the general ideal sheaf is equivalent to requiring that \[ (\mu-\chi) H^i(E(-2)) \to (\mu-2\chi) H^i(E(-1)) \] be an isomorphism. If we let $\cF'$ be the cokernel of a generic map \[ (\mu-\chi) \cO_{\bP^2}(-2) \to (\mu-2\chi) \cO_{\bP^2}(-1) \oplus \chi \cO_{\bP^2} \] Then the same argument as above shows that $H^i(\cF' \otimes E)=0$ if and only if the maps  \[(\mu-\chi) H^i(E(-2)) \to (\chi+1) H^i(E) \oplus (\mu-2\chi) H^i(E(-1))  \] are isomorphisms. But we've seen above that for general $E$ and a general map \[ (\mu-\chi) \cO_{\bP^2}(-2) \to (\mu-2\chi) \cO_{\bP^2}(-1) \oplus (\chi+1) \cO_{\bP^2} \] this is the case, and hence the same is true for a general map  \[ (\mu-\chi) \cO_{\bP^2}(-2) \to (\mu-2\chi) \cO_{\bP^2}(-1) \oplus \chi \cO_{\bP^2} \]

\end{proof}

Theorem 7.1 of \citep{huizenga} tells us when a general bundle of the above form has interpolation with respect to a general ideal sheaf of $n$ points.

\begin{thm}
Let $\alpha=1-\frac{\chi}{\mu}$. Then there is a bundle $E$ of the above form such that $E(\mu)$ has the interpolation property with respect to the Hilbert scheme of $n$ points in $\bP^2$ if and only if either $\alpha>\varphi^{-1}$, where $\varphi=\frac{1+\sqrt{5}}{2}$ or \[ \alpha \in \left \{ \frac{0}{1},\frac{1}{2},\frac{3}{5},\frac{8}{13},\ldots \right \} \] i.e.~$\alpha$ is a ratio of consecutive Fibonacci numbers.
\end{thm}

\section{Moving Curves}

In this previous section, we constructed effective divisors $D$ on \(N(\mu,\chi\) for certain choices of $\mu$ and $\chi$. In this section, we will show that these divisors are on the edge of the effective cone. For this purpose, we will construct a family of numerically equivalent curves on \(N(\mu,\chi)\) which pass through a general of $N(\mu,\chi)$ and have intersection number 0 with $D$. This will force $D$ to be on the boundary of the effective cone, but since the Picard number of $N(\mu,\chi)$ is two, this means that $D$ will generate an extremal ray of the effective cone. We then need only show that $D$ is not numerically proportional to $\cL_0$.

Take a general pencil of plane curves of degree $\mu$. $X$, the total space of the pencil, is isomorphic to the blowup of $\bP^2$ at the $\mu^2$ basepoints of the pencil. Given a line bundle $L$ on $X$ which has Euler characteristic $\chi$ when restricted to each fiber, we can push it forward to $\bP^2 \times \bP^1$ to get a family of semistable sheaves on $\bP^2$ parametrized by $\bP^1$, which in turn gives us a curve in $N(\mu,\chi)$.


\begin{lem}\label{vbintersect}

Suppose that $E$ is a vector bundle with the interpolation property with respect to $N(\mu,\chi)$. Let $D$ be the corresponding line bundle. Let $C$ be as above. Then we have \[ C \cdot (-D)=\chi(E)+\rk(E)\left(\frac{1}{2} \sum (a_i-a_i^2)+\frac{1}{2}b(b+3)\right)+b \frac{\ch_1(E)}{\rk(E)} \] 

\end{lem}
\begin{proof}
We again use Grothendieck-Riemann-Roch, but now we have \[\begin{array}{c}\ch(\pi_!(L \otimes E)) \todd(\bP^1)=\\\pi_*[(1+\sum a_i E_i+bH + \frac{1}{2}\left( b^2-\sum a_i^2 \right) H^2)(1+\frac{3}{2}H-\frac{1}{2}E+H^2)\ch(E)]\end{array} \] which equals \[ \chi(E)+\rk(E)\left(\frac{1}{2} \sum (a_i-a_i^2)+\frac{1}{2}b(b+3)\right)+b \frac{\ch_1(E)}{\rk(E)} \] Again, multiplying by the inverse of $\todd(\bP^1)$ does not change this, and since $D$ is given by the dual of the determinant line bundle corresponding to $E$, we have the desired result.
\end{proof}

\begin{lem}
Let $0<\chi \leq \mu$. Let \[ L=\sum_{i=1}^{\chi+\frac{1}{2}\mu(\mu-3)} E_i \] be a line bundle on $X$, and $C$ the corresponding curve on $N(\mu,\chi)$. Then curves in the numerical class of $C$ pass through a general point of $N(\mu,\chi)$.
\end{lem}
\begin{proof}
Let $M$ be a general line bundle of Euler characteristic $\chi$ on a general plane curve of degree $\mu$. Since $\chi(M)>0$, $M$ is linearly equivalent to a sum of (distinct) points $p_1,\ldots,p_n$. We need to show that there is a pencil of curves containing $p_1,\ldots,p_n$, but containing $p_i$ imposes one linear condition on $\bP^N$ where \[N={\frac{1}{2}(\mu+1)(\mu+2)-1}\] and \[\frac{1}{2}(\mu+1)(\mu+2)-1-\frac{1}{2} \mu(\mu-3)-\chi=3\mu-\chi\] will be bigger than 1 for $\chi \leq \mu$.
\end{proof}

Lemma \ref{vbintersect} lets us work out that the corresponding curves $C$ have intersection number 0 with the divisors associated to the vector bundles of theorem \ref{steiner}, since these have Euler characteristic 0 and $a_i=a_i^2$, but by the discussion at the beginning of the section, this means that the divisor $D$ is on the edge of the effective cone.

In the case of $\chi=0$, we must use a slightly different moving curve. Again, we take a general pencil of plane curves of degree $\mu$, but now we take the line bundle \[ -E_k+\sum_{i=1}^{\frac{1}{2}(\mu-1)(\mu-2)} E_i \] where $E_k$ is some exceptional divisor which is not included in the second part of the expression. The same arguments as above show that the divisor $D$ gives an edge of the effective cone. We can now prove the following theorem.

\begin{thm}
Let $D$ be the locus in $N(\mu,\chi)$ of sheaves which are not cohomologically orthogonal to a vector bundle coming from theorem \ref{steiner}, if such a vector bundle exists. Then $D$ and $\cL_0$ span the two edges of the effective cone of $N(\mu,\chi)$. More explicitly, the effective cone is spanned by $\cL_0$ and $\cL_1+\frac{\chi-\mu}{\delta} \cL_0$.
\end{thm}
\begin{proof}
We just need to show that $D$ is not proportional to $\cL_0$. But the moving curves we constructed above have intersection number 0 with $D$ and intersection number 1 with $\cL_0$, which would be impossible if the two were proportional.

The description of $D$ as a linear combination of $\cL_0$ and $\cL_1$ comes from the discussion following the statement of theorem \ref{steiner}.
\end{proof}

\section{Bridgeland Stability Conditions}

In this section, we will give a brief review of the theory of Bridgeland stability conditions, and prove one new result about them. I will restrict myself to the case of $\bP^2$, and I will often ignore the case of torsion-free sheaves. I will largely follow the treatment of \citep{abch} and \citep{bridgeland}.

For us, $D^b(\bP^2)$ will mean the bounded derived category of coherent sheaves on $\bP^2$. By $\cH^i$, we will mean the $i$th cohomology sheaf of an object of $D^b(\bP^2)$.

\begin{defn}
The heart of a bounded t-structure on $D^b(\bP^2)$ is a full additive subcategory $\cA \subset D^b(\bP^2)$ such that for all $A,B \in \cA$, $\Hom(A,B[k])=0$ if $k<0$ and for any object $E \in D^b(\bP^2)$ there are objects $0=E_m,E_{m+1},\ldots,E_n=E$ and triangles $E_i \to E_{i+1} \to F_i \to E_i[1]$ such that $F_i[i] \in \cA$.
\end{defn}

The standard example of the heart of a bounded t-structure on $D^b(\bP^2)$ is the full subcategory of coherent sheaves. In general, many of the things which make sense for the category of coherent sheaves work just as well for the heart of any bounded t-structure. For example, a map $A \to B$ with $A,B \in \cA$ is an inclusion with respect to $\cA$ if the mapping cone is also in $\cA$.

\begin{defn}
A stability condition on $\bP^2$ is a function $Z:K(\bP^2) \to \C$ called a central charge, together with $\cA$, the heart of a bounded t-structure on $D^b(\bP^2)$, which satisfy the following conditions. First, for any nonzero $A \in \cA$, $\Arg(Z(A)) \in \R_{>0} e^{i \theta}$ with $0<\theta \leq \pi$. We note that this allows us to put a partial order on the arguments which occur, with one argument being bigger than another if it is closer to the negative real axis.

We will call a nonzero $A \in \cA$ semistable with respect to the stability condition if every nonzero proper subobject $B \subset A$ has $\Arg(Z(B)) \leq \Arg(Z(A))$. We will call $A$ stable if this inequality is always strict.

The second part of being a stability condition is that any nonzero object $A \in \cA$ has a Harder-Narasimhan filtration, i.e.~a finite filtration \[ 0=\cF_0 \subset \cF_1 \subset \cdots \subset \cF_{n-1} \subset \cF_n=A \] such that for each of the subquotients $E_i=\cF_i/\cF_{i-1}$, $\Arg(Z(E_i))<\Arg(Z(E_{i+1}))$ and the $E_i$ are semistable.
\end{defn}

The Harder-Narasimhan filtration is always unique. Any semistable object has a Jordan-H\"older filtration, which is a finite filtration where the subquotients $E_i$ are stable and $\arg(Z(E_i))=\arg(Z(E_j))$ for all $i$ and $j$. The Jordan-H\"older filtration is not unique, but the subquotients are well-defined up to reordering. In general, we will say that two objects are S-equivalent with respect to some stability condition if the isomorphism classes of the subquotients in their Jordan-H\"older filtration are the same up to reordering. We note that if we have an inclusion of two semistable objects $E \subset F$ with $\arg(Z(E))=\arg(Z(F))$, then $E/F$ if semistable and $\arg(Z(E/F))=\arg(Z(E))$.

Given a stability condition, we can define a slope function $\mu_Z(E)=\frac{-\re(Z(E))}{\im(Z(E))}$. An object $A \in \cA$ will then be semistable if and only if for all nonzero proper subobjects $B \subset A$, $\mu_Z(B) \leq \mu_Z(A)$, and it will be stable if and only if this inequality is always strict.

There is another, equivalent, definition of stability condition, which will also prove useful.

\begin{defn}
A slicing of $D^b(\bP^2)$ is a collection of full additive subcategories $\cA_\phi$ with $\phi \in \R$, the objects of which will be called semistable of phase $\phi$, satisfying the following conditions.
\begin{enumerate}
\item For $A \in \cA_\phi$, $B \in \cA_{\phi'}$ with $\phi'<\phi$, $\Hom(B,A)=0$.
\item $\cA_\phi[1]=\cA_{\phi+1}$
\item For each nonzero object $C \in D^b(\bP^2)$, we have a collection of objects and maps \[ 0=E_0 \to E_1 \to \cdots \to E_n=C \] such that the mapping cone of $E_i \to E_{i+1}$ is in $\cA_{\phi_i}$ with the $\phi_i>\phi_{i+1}$. This collection of maps is called the Harder-Narasimhan filtration of $C$ with respect to the slicing. These mapping cones will be called the subquotients of the filtration.
\end{enumerate}
\end{defn}

The following is proposition 5.3 of \citep{bridgeland}.

\begin{prop}
A stability condition on $D^b(\bP^2)$ is equivalent to a slicing, together with a central charge function $Z:K(\bP^2) \to \C$ such that for any nonzero $A \in \cA_\phi$, $\arg(Z(A))=\pi \phi$.
\end{prop}

The idea behind the correspondence is as follows. Given a slicing of $D^b(\bP^2)$, we can take the full extension-closed subcategory of $D^b(\bP^2)$ generated by the objects of $\cA_\phi$ with $\phi \in (0,1]$. This will be the heart of a bounded t-structure. To go the other way, for $\phi \in (0,1]$, we let $\cA_\phi$ be the full subcategory of objects $E$ in the heart which are semistable and such that $\arg(Z(E))=\pi \phi$.

\begin{defn}
Given a slicing, for $I \subset \R$ an interval, we can define $\cA_I$ to be the extension-closed subcategory of $D^b(\bP^2)$ generated by objects in $\cA_\phi$ with $\phi \in I$.

A slicing is called locally finite if for all $\phi$, there is some $\epsilon>0$ such that $\cA_{(\phi-\epsilon,\phi+\epsilon)}$ is of finite length. A stability condition is called locally finite if the corresponding slicing is locally finite.
\end{defn}

In general, we can define a metric on the set of (locally finite) stability conditions as follows.

\begin{defn}
Suppose we have two slicings $\cP$ and $\cP'$ of $D^b(\bP^2)$. Given a nonzero object $E \in D^b(\bP^2)$, we will define $\phi^+_{\cP}(E)$ (resp.~$\phi^-_{\cP(E)}$) to be the subquotient of the $\cP$-Harder-Narasimhan filtration of $E$ with largest (resp.~smallest) phase. We define \[ d(\cP,\cP')=\sup_{0 \neq E \in D^b(\bP^2)} \max\left \{|\phi_{\cP}^+(E)-\phi_{\cP'}^+(E)|,|\phi_{\cP}^-(E)-\phi_{\cP'}^-(E)| \right \} \]

Given two stability conditions $(Z,\cP)$ and $(Z',\cP')$, we define $d((Z,\cP),(Z',\cP'))=d(\cP,\cP')+||Z-Z'||$, where for the last term of the sum, we use the Euclidean norm on the finite-dimensional vector space $\Hom(K(\bP^2),\C)$.
\end{defn}

Theorem 7.1 of \citep{bridgeland} says that the map from the space of (locally finite) stability conditions with the above metric to $\Hom(K(\bP^2),\C)$ is a local homeomorphism onto a linear subspace.

On $\bP^2$, there is a family of (locally finite) stability conditions parametrized by two variables $s \in \R$ and $t \in \R_{>0}$, which we will think of as giving coordinates on the upper half plane. The central charge is given by \[ Z_{s,t}(E)=-\int_{\bP^2} e^{-(s+it)H} \chern(E) \] 

The heart of the t-structure is a little more complicated. Given a nonzero torsion-free sheaf $E$ on $\bP^2$, we will define its slope $\mu(E)$ to be $\frac{c_1(E)H}{\rk(E)}$. We will say that $E$ is $\mu$-semistable if for all nonzero proper subobjects $F \subset E$, we have $\mu(F) \leq \mu(E)$. This lets us define a notion of Harder-Narasimhan filtration, which we will call the $\mu$-Harder-Narasimhan filtration. 

\begin{defn}
For $s \in \R$, define $\cF_s$ to be the full subcategory of coherent sheaves $A$ on $\bP^2$ which are torsion-free and such that for each subquotient $A_i$ of $A$ in the $\mu$-Harder-Narasimhan filtration of $A$, $\mu(A_i) \leq s$.

We will define $\cQ_s$ to be the full subcategory of coherent sheaves $A$ on $\bP^2$ such that for each subquotient $A_i$ of $A/\tor(A)$ (where $\tor(A)$ denotes the torsion subsheaf of $A$) in the $\mu$-Harder-Narasimhan filtration of $A/\tor(A)$, $\mu(A_i)>s$.

We define a category $\cA_s \subset D^b(\bP^2)$ as the full subcategory containing objects $C$ with $\cH^{-1}(C) \in \cF_s$, $\cH^0(C) \in \cQ_s$, and $\cH^i(C)=0$ for $i \neq 0,1$.
\end{defn}

The following is implied by theorem 5.11 of \citep{abch} and the preceding discussion.

\begin{thm}
For all $s \in \R$ and $t>0$, the category $\cA_s$ is the heart of a bounded t-structure on $D^b(\bP^2)$ and the pair $(Z_{s,t},\cA_s)$ give a stability condition.
\end{thm}

We note that this gives a set-theoretic inclusion from the upper half plane to the space of stability conditions. By inspection, the composite map to $\Hom(K(\bP^2),\C)$ is a local embedding, so the map from the upper half plane to the space of stability conditions must be continuous.

It will be useful to have a more explicit description of the slope function at each point of the upper half plane. We have \[ \mu_{s,t}(A)=\frac{\frac{\rk(A)}{2}(s^2-t^2)+\ch_2(A)-s \ch_1(A)}{t(\ch_1(A)-s \rk(A))} \]

Given a sheaf $\cF$ (or more generally, an element of $K(\bP^2)$), we get a collection of potential walls in the upper half plane. Each potential wall is the locus where $\mu_{s,t}(\cF)=\mu_{s,t}(E)$ for some fixed sheaf $E$. If $\cF \in N(\mu,\chi)$, then these potential walls are concentric semicircles. More precisely, we have the following result of \citep{abch}.

\begin{lem}
The potential wall corresponding to a sheaf $E$ is a semicircle with center \[ \left(\frac{\chern_2(\cF)}{\chern_1(\cF)},0\right)=\left(\frac{\chi}{\mu}-\frac{3}{2},0\right) \] and radius \[ \sqrt{\left (\frac{\chi}{\mu}-\frac{3}{2} \right )^2+\frac{2}{\chern_0(E)} \left ( \chern_2(E)-\left (\frac{\chi}{\mu}-\frac{3}{2}\right)\chern_1(E) \right)} \] In particular, the potential walls never intersect.
\end{lem}

The following proposition allows us to use the geometry of the potential walls to understand when sheaves have nonzero maps to a torsion sheaf. It extends lemma 6.3 of \citep{abch} to the case of torsion sheaves.

\begin{prop}\label{wallconsistency}
Suppose we have a map $E \to \cF$, with $\cF \in N(\mu,\chi)$. Suppose at some point $(s,t)$ on the potential wall corresponding to $E$, $E \to \cF$ is an inclusion of semistable objects in the heart of the t-structure. Then this is true for all $(s',t')$ on that potential wall.
\end{prop}
\begin{proof}

Let $C$ be the mapping cone of $E \to \cF$. We have a distinguished triangle \[ E \to \cF \to C \to E[1] \] Since $\cF$ is a sheaf, and both $E$ and $C$ are in $\cA_s$, the corresponding long exact sequence of cohomology sheaves shows that $E$ is a sheaf too. By assumption, with respect to $(s,t)$; $E$, $\cF$, and $C$ are all in the same slice, say $\cA_\phi$. For pure one-dimensional sheaves and $t>0$, we have $\phi \in (0,1)$.

In order to prove the theorem, it is enough to show that for all $(s',t')$ on the wall; $E$, $\cF$, and $C$ will all belong to $\cA'_\phi$, that is, the corresponding slice with respect to the slicing given by $(s',t')$. Since $\phi \in (0,1]$, this will imply that $E \to \cF$ is an inclusion of objects in the heart, and the fact that they're in the corresponding slice implies that they are all semistable of the same slope.

In order to show this, it is enough to show that $E$, $\cF$, and $C$ all belong to $\cA'_\phi$ for all stability condition $(s',t')$ on the wall which are contained in a ball of radius $\frac{\epsilon}{2}$ centered at $(s,t)$, where the ball is defined with respect to the above metric on the space of stability conditions and $\epsilon<\min\{phi,1-\phi\}$. This is because the metric on the space of stability conditions defines the usual topology on the wall, so any two points $(s,t)$ and $(s',t')$ on the wall can be connected by a sequence of points $(s,t)=(s_0,t_0),(s_1,t_1),\ldots,(s_n,t_n)=(s',t')$ such that $d((s_i,t_i),(s_{i+1},t_{i+1}))<\frac{\epsilon}{2}$.

Suppose that there are two stability conditions $(s,t)$ and $(s',t')$ within $\frac{\phi}{2}$ of each other such that $E$, $\cF$, and $C$ are all in the same $(s,t)$-slice but not the same $(s',t')$-slice. By definition of the potential wall, the argument of the central charge will be the same, so $E$, $\cF$, or $C$ must either be semistable in $\cA_{\phi+n}$ for some $n \in 2\Z$, or one of them must stop being semistable. By definition of the metric on the stability manifold, for $(s',t')$ within $\frac{\phi}{2}$ of $(s,t)$, the first case is impossible.

We claim that if $\cF$ remains semistable, then so must $E$ and $C$. Suppose $E$ stops being semistable. Let $E'$ be the semistable factor with the largest slice in the Harder-Narasimhan filtration with respect to the $(s',t')$-slicing. Then the phase of $E'$ must be larger than $\phi$, since the overall phase of $E$ is $\phi$. By definition of a slicing, this means that the map $E' \to \cF$ is 0. By applying $\Hom(E',\cdot)$ to the above triangle, we see that there must be a nonzero map $E' \to C[-1]$. With respect to $(s,t)$, $C[-1] \in \cA_{\phi-1}$, and for $(s',t')$ within $\frac{\epsilon}{2}$ of $(s,t)$, the largest subquotient occurring in the $(s',t')$-Harder-Narasimhan filtration of $C[-1]$ will have a $(s',t')$-phase less than $\phi-\frac{1}{2}$, but this means that any map $E' \to C[-1]$ vanishes, a contradiction. A similar argument shows that if $\cF$ continues to be semistable, then so must $C$.

We now want to show that $\cF$ is $(s',t')$-semistable of phase $\phi$. Suppose it's not. Let $\cF'$ be the piece of largest phase in the $(s',t')$-Harder-Narasimhan filtration of $\cF$. By nestedness of potential walls, the $(s,t)$-slope of $\cF'$ is still greater than the $(s,t)$-slope of $\cF$. Let $\cG$ be piece of the $(s,t)$-Harder-Narasimhan filtration of $\cF'$ with greatest $(s,t)$-phase. By our choice of $\epsilon$, $\cG$ will be contained in $(0,1]$. Since the $(s,t)$-slope of $\cG$ is bigger than the $(s,t)$-slope of $\cF'$, we can conclude that the $(s,t)$-phase of $\cG$ is bigger than $\phi$.

Since $\cF$ is $(s,t)$-semistable of phase $\phi$, we must have the map $\cG \to \cF$ is zero. Let $K$ be the mapping cone of $\cF' \to \cF$. Then we have a distinguished triangle \[ K[-1] \to \cF' \to \cF \to K \]

Applying $\Hom(\cG,\cdot)$ to this triangle, and using the fact that $\Hom(\cG,\cF)=0$, but $\Hom(\cG,\cF') \neq 0$, we see that there must be a nonzero map $\cG \to K[-1]$. We know that the largest piece of the $(s',t')$-Harder-Narasimhan filtration of $K$ must have phase at most $\phi+\frac{\epsilon}{2}$, so the largest piece of the $(s,t)$-Harder-Narasimhan filtration of $K[-1]$ can have phase at most $\phi-\frac{1}{2}$, and so there can be no nonzero maps $\cG \to K[-1]$, a contradiction.
\end{proof}

\section{Nef Cone}

We recall the following from \citep{bayermacri}.

\begin{thm}
Let $Z$ be a Bridgeland stability condition on $D^b(\bP^2)$. Let $v \in K(\bP^2)$. Let $C$ be a curve, and $\sC \in D^b(C \times \bP^2)$. We can associate to this a real number \[\Im \left ( \frac{-Z(\phi_{\sC}(\cO_C))}{Z(v)} \right ) \] where $\phi_{\sC}:D^b(C) \to D^b(\bP^2)$ denotes the integral transform with kernel $\sC$.

This assignment induces a map $N_1(\cM_Z(v)) \to \R$, or equivalently, gives an element of $N^1(\cM_Z(v))$. This element is nef, and it has intersection number 0 with a curve if and only if the curve parametrizes a family of strictly $Z$-semistable complexes all of which are $S$-equivalent.
\end{thm}

We already know that $\cL_0$ is on the edge of the nef cone, so it remains to find the other edge of the nef cone. To do this, we will find a Bridgeland stability condition $(s,t)$ such that all Simpson-semistable sheaves are $(s,t)$-semistable, and there is a one-dimensional family of sheaves which are $S$-equivalent with respect to the $(s,t)$-stability condition, but not with respect to Simpson stability. The above theorem will then give us a divisor class which is nef but has intersection number 0 with some curve, and hence must lie on the edge of the nef cone. It will then suffice to check that this divisor class is not a multiple of $\cL_0$. By theorem \ref{moridreamspace}, this divisor class will be a positive real multiple of the class of a semiample line bundle.

\begin{lem}
Let $\cF$, $\cG$ be pure one-dimensional sheaves, and $(s,t)$ a stability condition. If $\cF \to \cG$ destabilizes $\cG$ with respect to $(s,t)$, then it also destabilizes $\cG$ with respect to Simpson stability.
\end{lem}

\begin{proof}
We have $\mu_{s,t}(\cF)=\frac{1}{t} \left ( \frac{\ch_2(\cF)}{\ch_1(\cF)-s} \right )$, and similarly for $\cG$. But $\ch_2(\cF)=\chi(\cF)-\frac{3}{2} \ch_1(\cF)$. This agrees with our earlier definition of Simpson stability.
\end{proof}

The long exact sequence of cohomology sheaves associated to a destabilizing triangle, along with the fact that all the objects in the triangle belong to the heart $\cA_s$ together imply that any Bridgeland-destabilizing object must itself be a sheaf. The next lemma shows that for sufficiently large $t$, torsion-freesheaves cannot be Bridgeland destabilizing. Combined with the previous lemma, this will show that for sufficiently large $t$, Simpson-stability implies Bridgeland stability, and similarly for semistability.

\begin{lem}
Fix the numerical invariants for a pure one-dimensional sheaf. Then the radius of the wall in the stability manifold given by a torsion-free sheaf $E$ is bounded. More specifically, the radius of the wall associated to $E$ is at most \[ \frac{\ch_1(\cF)}{2 \rk(E)}\]
\end{lem}
\begin{proof}

Let $d=\ch_2(\cF)=\chi-\frac{3}{2}\mu$. Suppose we also fixed the Chern character of the destabilizing object $E$ to be $(r',c',d')$. Then the potential wall is a semicircle with center \[ \left ( \frac{d}{\
mu},0 \right ) \] and radius \[ \sqrt{\frac{d^2}{\mu^2}+\frac{2d'}{r'}-\frac{2c'd}{r'\mu}} \] If this is an actual wall, then by proposition \ref{wallconsistency}, $E \to \cF$ is an inclusion of objects in $\cA_s$ at every point of the wall.

Let us denote by $R$ the radius of the wall, and $K$ the kernel of the map $E \to \cF$. We have the string of inequalities \[ \frac{d}{\mu}+R < \frac{c_1(E)}{\rk(E)} \leq \frac{\mu}{\rk(E)}+ \frac{c_1(K)}{\rk(K)} \leq \frac{\mu}{\rk(E)}+\frac{d}{\mu}-R \] which together give us the inequality \[ R< \frac{\mu}{2 \rk(E)} \]
\end{proof}

The following lemma follows from the calculation in \citep{abch} of the nef cone of the Hilbert scheme together with the correspondence between chambers in the space of stability conditions on $\bP^2$ and cones in the Picard group of the Hilbert scheme.

\begin{lem}
Let $\cI_W(k)$ be a twist of an ideal sheaf. Then $\cI_W(k)$ is $(s,t)$-semistable if $(s,t)$ is to the left of the line $s=k$ and $(s,t)$ is outside the potential wall of $\cO(k-1)$, which is a semicircle with center \[ \left (k-|W|-\frac{1}{2},0\right )\] and radius $|W|-\frac{1}{2}$
\end{lem}

We will now find the destabilizing object which gives the largest wall for sheaves in $N(\mu,\chi)$. Write $a=d+\frac{1}{2} \mu^2$. If $\cF$ is the pushforward of a line bundle of degree $a$ on a smooth plane curve of degree $\mu$, then $\cF \in N(\mu,\chi)$. Now write $a=b\mu+\epsilon$ with $-\frac{\mu}{2} \leq \epsilon \leq \frac{\mu}{2}$. Note that the case $\epsilon=\pm \frac{\mu}{2}$ is ambiguous. Since we will be looking at many rank one walls, we note that the potential wall associated to $\cI_W(k)$ is a semicircle with center \[ \left ( b+\frac{\epsilon}{\mu}-\frac{\mu}{2},0 \right ) \] and radius \[ \sqrt{\left (\frac{d}{\mu}-k \right)^2-2|W|}=\sqrt{\left(\frac{\epsilon}{\mu}-\frac{\mu}{2}+(b-k)\right)^2-2|W|} \]

\begin{prop}
If $\epsilon \leq 0$, then $\cI_W(b)$ is destabilizing and gives the biggest wall, where $|W|=|\epsilon|$. If $\epsilon \geq 0$, then $\cO(b)$ is destabilizing and gives the biggest wall.
\end{prop}

\begin{proof}

In order to prove this, there are a number of things we must check. First, we must check that for some one-dimensional sheaf $\cF$, this is an inclusion of semistable objects of the same phase. We first consider the case $\epsilon \leq 0$. We want to check that on the potential wall corresponding to $\cI_W(b)$, the ideal sheaves are semistable, and there is an inclusion of objects in $\cA_s$.

The potential wall associated to $\cI_W(b)$ with $|W|=-\epsilon$ has radius $\frac{\epsilon}{\mu}+\frac{\mu}{2}$. We know that the points $\left ( b+\frac{2\epsilon}{\mu},0 \right )$ and $(b-\mu,0)$ are on the potential wall. These are to the left of the line $s=b$. They strictly contain the semicircle with center $(b-\epsilon-\frac{1}{2},0)$ and radius $\epsilon-\frac{1}{2}$ unless $\epsilon=-\frac{\mu}{2}$, in which case the semicircles are equal.

Let $C$ be a smooth plane curve. Consider the line bundle $bH-p_1 \cdots -p_\epsilon$. This gives a semistable sheaf $\cF$ of the right numerical invariants, and we have a map $\cI_{p_1,\ldots,p_\epsilon}(b) \to \cF$ which is surjective with kernel  $\cO(b-\mu)$. Since the left edge of the wall is $(b-\mu,0)$, we see that when $t>0$, the mapping cone will be in $\cA_s$.

The case where $\epsilon \geq 0$ is very similar, but slightly easier. In this case, the destabilizing object is the line bundle $\cO(b)$, which is stable if $s<b$. The potential wall associated to $\cO(b)$ has radius $\frac{\mu}{2}-\frac{\epsilon}{\mu}$, so we see that $\cO(b)$ will be stable along this potential wall, since the right edge is $(b,0)$. If we have a line bundle of the form $bH+p_1+\cdots+p_\epsilon$, then we get a map $\cO(b) \to \cF$ which has kernel $\cO(b-\mu)$ and cokernel a zero-dimensional sheaf. Since the left edge of the potential wall is $(b-\mu+\frac{2\epsilon}{\mu},0)$, this will always be in the heart at every point of the potential wall.

We must now check that there can be no bigger destabilizing walls. This will also imply that $\cI_W(b)$ really destabilizes the sheaves described in the previous paragraph. The largest possible radius of a higher rank wall is $\frac{\mu}{4}$, but when $\mu>2$, we have $\frac{\epsilon}{\mu}+\frac{\mu}{2} > \frac{\mu}{4}$, so all higher rank walls are strictly smaller.

Now suppose that $\cI_Z(k)$ gives a bigger wall. Since this must be in $\cA_s$ for $s=b+\frac{2\epsilon}{\mu}$, we must have $k \geq b+\frac{2\epsilon}{\mu}$. Let $K$ be the kernel of the map $\cI_Z(k) \to \cF$. We have that $K$ is a torsion-free rank 1 sheaf, and $k-\mu \leq c_1(K) \leq k$. Since the mapping cone of $\cI_Z(k) \to \cF$ must be in $\cA_s$ when $s=b-\mu$, we have that $c_1(K) \leq b-\mu$, and so $k \leq b$. We see that $b=k$ unless $\epsilon=\frac{-\mu}{2}$, in which case $b$ can be $k-1$. We note that when $\epsilon=-\frac{\mu}{2}$, both $\cO(b-1)$ and $\cI_W(b)$ with $|W|=\frac{\mu}{2}$ give the same semicircle.

Increasing the number of points in the ideal sheaf decreases the radius of the potential wall, so it suffices to show that we can't have $\cI_Z(b)$ with $|Z|<|W|$. In the case of $\epsilon \geq 0$ (and hence also $\epsilon=-\frac{\mu}{2})$, the destabilizing object is a line bundle, so this case is complete. From the inequalities above, we know that $c_1(K) \leq b-\mu$. But if $\epsilon \leq 0$, then $(b-\mu,0)$ is on the left edge of the wall associated to $\cI_W(b)$, so for a bigger wall, the mapping cone will not be in the heart along the leftmost part of the wall.

\end{proof}

The above proof also gives the following result.

\begin{cor}
Suppose $-\frac{\mu}{2}<\frac{\epsilon}{\mu} \leq 0$. Then $\cI_W(b)$ with $|W|=|\epsilon|$ are the only destabilizing objects along the largest wall. The Jordan-H\"older factors of $\cF$ for sheaves which become semistable along this wall are $\cI_W(b)$ and $\cO(b-\mu)[1]$. If $0 \leq \frac{\epsilon}{\mu} < \frac{\mu}{2}$, then the only destabilizing objects along the largest wall are of the form $\cO(b)$, and the Jordan-H\"older subquotients are $\cO(b)$ and an extension of a zero-dimensional sheaf by $\cO(b-c)[1]$. If $\frac{\epsilon}{\mu}=\pm \frac{\mu}{2}$, then the destabilizing objects can be can be of either form. In either case, the Jordan-H\"older subquotients are again of the same form as above, unless all the points in the ideal sheaf or the zero-dimensional sheaf are collinear, in which case the Jordan-H\"older quotients are $\cO(b-1)$, $\cO(b-\mu)[1]$, and $\cO_\ell(b-\frac{\mu}{2})$, where $\ell$ is the line containing all the points.
\end{cor}

From the above, we see that the pushforward of a line bundle on a plane curve becomes semistable if the difference with the closest multiple of the hyperplane class is effective. In particular, the general one-dimensional sheaf does not become semistable unless $\mu=3$, because this is the only case when a general line bundle of degree 1 is effective. In particular, this line bundle defines a birational morphism, so it cannot be a multiple of $\cL_0$. In the case of $c=3$, we see that we get a morphism $N(3,1) \to \bP^2$, which for a smooth plane curve remembers the point defining the corresponding line bundle on that curve.

\section{Isomorphisms of Moduli Spaces of One-Dimensional Sheaves}

When are two moduli spaces of one-dimensional sheaves $N(\mu,\chi)$ and $N(\mu',\chi')$ isomorphic? The dimension of $N(\mu,\chi)$ is $\mu^2+1$, so a necessary condition is that $\mu=\mu'$. If $\mu=\mu'$ and $\chi \cong \pm \chi' \pmod \mu$, then we have seen that $N(\mu,\chi) \cong N(\mu,\chi')$ in corollary \ref{isomorphisms}.

This condition is always satisfied for $\mu=1$. For $\mu=2$, we have $N(\mu,0) \cong N(\mu,1) \cong \bP^5$, so this condition is not necessary. However, we will show that this is the only case.

\begin{thm}
Let $\mu \geq 3$. Then $N(\mu,\chi) \cong N(\mu,\chi')$ only if $\chi \cong \pm \chi' \pmod \mu$.
\end{thm}
\begin{proof}

Let $a$ be as in the previous section, and $a'$ the analogue for $\chi'$. Then $a=\frac{1}{2}\mu(\mu-3)+\chi$, so $\chi \cong \pm \chi' \pmod \mu$ if and only if $a \cong \pm a' \pmod \mu$. Write $a=b\mu+\epsilon$ with $\frac{-\mu}{2} \leq \epsilon \leq \frac{\mu}{2}$. We note that there is a possible ambiguity if $\epsilon = \frac{\mu}{2}$. Similarly, write $a'=b' \mu+\epsilon'$. It suffices to show that we must have $|\epsilon|=|\epsilon'|$.

Let us suppose that $N(\mu,\chi) \cong N(\mu,\chi')$. This isomorphism induces an isomorphism of Picard groups which must preserve the nef cone. Since one edge of the nef cone is distinguished by the fact that it is proportional to the canonical class, this isomorphism must preserve the other edge of the nef cone too.

This implies that both the map $s$ which sends a sheaf to it support and the map $f$ given by the other edge of the nef cone are preserved by isomorphisms, as is $E$, be the exceptional locus of $f$. Let us consider the restriction of $s$ to $E$. By the discussion in the above section, the general fiber has dimension $|\epsilon|$, since we can identify the fiber over a point in $\bP H^0(\cO_{\bP^2}(\mu)))$ corresponding to a smooth plane curve $C$ with the $|\epsilon|$-th symmetric power of $C$ as embedded in its Jacobian by the Abel-Jacobi map. This shows that $|\epsilon|$ must equal $|\epsilon'|$
\end{proof}



\bibliography{torsionbib}

\end{document}